\documentclass[11pt]{article}
\usepackage{mathrsfs}
\usepackage[centertags]{amsmath}
\usepackage{amsfonts}
\usepackage{amssymb}
\usepackage{amsthm}
\usepackage{cases}
\usepackage{indentfirst}
\usepackage{epsfig}
\usepackage{graphicx}
\usepackage{color}
\usepackage{CJK}
\usepackage[numbers,sort&compress]{natbib}
\usepackage{subfigure} 

\date{}

\newtheorem{theorem}{Theorem}[section]

\newtheorem{lemma}{Lemma}[section]

\newtheorem{proposition}{Proposition}[section]
\newtheorem{remark}{Remark}[section]

\numberwithin{equation}{section}

\allowdisplaybreaks[4]

\begin{document}
\title{\textbf{On the distributional expansions of powered extremes from Maxwell distribution}\thanks{
{Corresponding author. E-mail address: wjj@swu.edu.cn (J. Wang).}
}}
\author{Jianwen Huang$^a$,\quad Xinling Liu$^b$, \quad Jianjun Wang$^a$,\quad Zhongquan Tan$^c$,\quad Jingyao Hou$^a$,\quad Hao Pu$^d$\\
{\small $^a$School of Mathematics and Statistics, Southwest
University, Chongqing, 400715, China}\\{\small $^b$School of Mathematics and Information, China West Normal
University, Nanchong, 637002, China}\\{\small $^c$Department of Statistics, Jiaxing University, Jiaxing, 314000, China}\\{\small $^d$School of
Mathematics, Zunyi Normal College, Zunyi,
563002, China}
}

\maketitle
\begin{quote}
{\bf Abstract.}~~In this paper, asymptotic expansions of the distributions and densities of powered extremes for Maxwell samples are considered. The results show that the convergence speeds of normalized partial maxima relies on the powered index. Additionally, compared with previous result, the convergence rate of the distribution of powered extreme from Maxwell samples is faster than that of its extreme. Finally, numerical analysis is conducted to illustrate our findings.

\end{quote}


\section{Introduction}
\label{sec1}

In extreme value theory, researchers recently focus on investigating the quality of convergence of normalized $\max\{X_k,1\leq k\leq n\}:=M_n$ of a sample. For the convergence rate of normalized $M_n$, general cases were discussed by Smith \cite{Smith}, Leadbetter et al. \cite{Leadbetter et al.}, Galambos \cite{Galambos} and de Haan and Resnick \cite{de Haan and Resnick}, and specific cases were considered by Hall \cite{Hall 1979,Hall 1980}, Nair \cite{Nair}, Liao and Peng \cite{Liao and Peng}, Lin et al. \cite{Lin et al 2011,Lin et al 2016}, Du and Chen \cite{Du and Chen,Chen and Du}, and Huang et al. \cite{Huang et al 2017}. Hall \cite{Hall 1980} derived the asymptotics of distribution of normalized $|M_n|^t$, the powered extremes for given power index $t>0$. Zhou and Ling \cite{Zhou and Ling} improved Hall' results and proved that the convergence speed of distributions and densities of extremes depends on the power index. Nair \cite{Nair} established the asymptotic expansions of normalized maximum from normal samples. Liao et al. \cite{Liao et al} and Jia et al. \cite{Jia et al} generalized Nair's work to skew-normal distribution and general error distribution, respectively.

Since the Maxwell distribution was proposed by James Clerk Maxwell \cite{Mandl}, a variety of applications of it in physics (in particular in statistical mechanics) have been found; see Shim and Gatignol \cite{Shim and Gatignol}, Tomer and Panwar \cite{Tomer and Panwar} and Shim \cite{Shim} and some statisticians and reliability engineers have investigated the statistical properties of it as well, see \cite{Liu and Liu,Huang and Chen,Dar et al,Huang et al 2017,Huang et al 2017b,Huang et al 2018,Huang and Wang 2018a,Huang and Wang 2018b}.

The aim of this paper is to investigate the distributional tail representation of $|X|^t$ with $X$ following Maxwell distribution and the limiting distribution of normalized $|M_n|^t$, and obtain asymptotic expansions of distribution and density of powered maximum from Maxwell distribution.

\begin{figure}[h]
\begin{center}
{\includegraphics[width=0.40\textwidth,height=40mm]{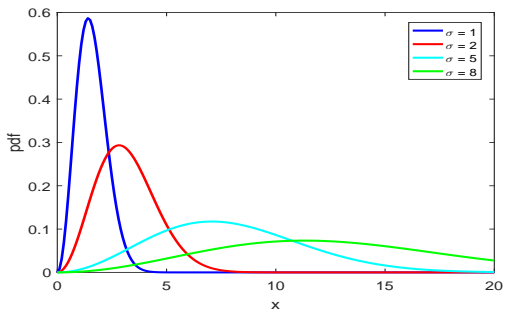}}
\caption{Probability density function of Maxwell distribution}\label{fig.1}
\end{center}
\vspace*{-14pt}
\end{figure}

Let $\{X_n,n\geq1\}$ be a sequence of independent identically distributed (i.i.d.) random variables with marginal cumulative distribution function (cdf) $F$ obeying the Maxwell distribution (abbreviated as $F\sim MD$), and as before let $M_n=\max\{X_i,1\leq i\leq n\}$ denote the partial maximum of $\{X_n,n\geq1\}$. The probability density function (pdf) of the MD is defined by
\begin{align}\label{eq.1}
f(x)=\sqrt{\frac{2}{\pi}}\frac{x^2}{\sigma^3}\exp\left(-\frac{x^2}{2\sigma^2}\right),~x>0,
\end{align}
where $\sigma>0$ is the scale parameter. Figure \ref{fig.1} presents the graph of pdf of Maxwell distribution. It shows that with the scale parameter increasing, the tail of pdf of MD becomes much heavier.

Liu and Liu \cite{Liu and Liu} showed that $F\in D(\Lambda)$, i.e., the max-domain of attraction of Gumbel extreme value distribution and the normalizing constants $a_n$ and $b_n$ can be given by
\begin{align}\label{eq.3}
a_n=\sigma^2b_n^{-1}
\end{align}
and
\begin{align}\label{eq.4}
\sqrt{\frac{\pi}{2}}\frac{\sigma}{b_n}\exp\left(\frac{b^2_n}{2\sigma^2}\right)=n
\end{align}
such that
\begin{align}\label{eq.5}
\lim_{n\to\infty}\mathbb{P}(M_n\leq a_nx+b_n)=\Lambda(x)=\exp\{-\exp(-x)\}.
\end{align}

The paper is constructed as follows. Section $2$ presents auxiliary lemmas with proofs. The main results are given in Section $3$. Numerical studies presented in Section $4$ compare the precision of the true values with its approximations. Section $5$ provides the proofs of main results.

\section{Auxiliary results}
\label{sec2}
To prove the main results, the following auxiliary lemmas are needed.
\begin{lemma}\label{le.1}
Let $F(x)$ and $f(x)$ respectively represent the cdf and the pdf of MD with $\sigma>0$, respectively. For large $x$, we have
\begin{align}\label{eq.2}
1-F(x)=\sigma^2x^{-1}f(x)\left[1+\sigma^2x^{-2}-\sigma^4x^{-4}+3\sigma^6x^{-6}+O(x^{-8})\right].
\end{align}
\end{lemma}
The proof of Lemma $\ref{le.1}$ is derived by integration by parts.

The following lemma gives the distributional tail representation of $X^t$ with $X\sim MD$.
\begin{lemma}\label{le.2}
Suppose that $0<t\neq2$. Let $F_t(x)$ denote the cdf of $X^t$ with $X\sim MD$. Then for large $x$, we get
\begin{align}\label{eq.6}
1-F_t(x)=C_t(x)\exp\left\{-\int^x_1\frac{g_t(u)}{\tilde{f}_t(u)}du\right\},
\end{align}
where
\begin{align}\notag
C_t(x)&\to\frac{2}{\sigma}\sqrt{\frac{2}{\pi}}\exp\left(-\frac{1}{2\sigma^2}\right)~\mbox{as}~x\to\infty,\\
\notag g_t(x)&=1-\sigma^2x^{-2/t}\to 1~\mbox{as}~x\to\infty,
\end{align}
and
\begin{align}\label{eq.7}
\tilde{f}_t(x)=\sigma^2tx^{1-\frac{2}{t}}~\mbox{with}~\tilde{f}_t'(x)\to 0~\mbox{as}~x\to\infty.
\end{align}
\end{lemma}

\noindent \textbf{Proof.} Combining with (\ref{eq.2}), we get
\begin{align}
\label{eq.8} 1-F_t(x)&=2\frac{\sigma^2f(x^{\frac{1}{t}})}{x^{\frac{1}{t}}}\left[1+\sigma^2x^{-\frac{2}{t}}
-\sigma^4x^{-\frac{4}{t}}+3\sigma^6x^{-\frac{6}{t}}+O(x^{-\frac{8}{t}})\right]\\
\notag&=\frac{2}{\sigma}\sqrt{\frac{2}{\pi}}\exp\left(-\frac{x^{\frac{2}{t}}}{2\sigma^2}+\frac{1}{t}\log x\right)\left[1+\sigma^2x^{-\frac{2}{t}}
-\sigma^4x^{-\frac{4}{t}}+3\sigma^6x^{-\frac{6}{t}}+O(x^{-\frac{8}{t}})\right]\\
\notag&=C_t(x)\exp\left(-\int^x_1\frac{g_t(u)}{\tilde{f}_t(u)}du\right)\left[1+\sigma^2x^{-\frac{2}{t}}
-\sigma^4x^{-\frac{4}{t}}+3\sigma^6x^{-\frac{6}{t}}+O(x^{-\frac{8}{t}})\right]
\end{align}
with $\tilde{f}_t(x)=\sigma^2tx^{1-\frac{2}{t}}$, $g_t(x)=1-\sigma^2x^{-2/t}$ and $C_t(x)\to\frac{2}{\sigma}\sqrt{\frac{2}{\pi}}\exp\left(-\frac{1}{2\sigma^2}\right)$ as $x\to\infty$.

\qed

Applying the result of Lemma $\ref{le.2}$ and Corollary $1.7$ \cite{Resnick}, the following result holds.
\begin{proposition}\label{pr.1}
Under the conditions of Lemma $\ref{le.2}$, we have $F_t(x)\in D(\Lambda)$, where $D(\Lambda)$ is the domain of $\Lambda(x)=\exp\{-\exp(-x)\}.$
\end{proposition}

Then, our aim is to select the suitable normalizing constants which ensure that the distribution of maximum tends to its extreme value limit. A combination of (\ref{eq.4}) and (\ref{eq.8}), we obtain that $d_n=b^t_n$. It follows from (\ref{eq.7}) that
\begin{align}\label{eq.12}
c_n=\tilde{f}_t(d_n)=\sigma^2tb^{t(1-\frac{2}{t})}_n=\sigma^2tb^{t-2}_n.
\end{align}

The following work is to find the special normalizing constants $c_n$ and $d_n$ for the case of powered index $t=2$. Similarly, it is necessary to establish the distributional tail representation of $X^2$ with $X\sim MD$.

\begin{lemma}\label{le.3}
Assume that $t=2$. Let $F_2(x)$ stand for the cdf of $X^2$ with $X\sim MD$. Then for large $x$, we get
\begin{align}\label{eq.9}
1-F_2(x)=C_2(x)\exp\left\{-\int^x_1\frac{g_2(u)}{\tilde{f}_2(u)}du\right\},
\end{align}
where
\begin{align}\notag
C_2(x)&\to\frac{2}{\sigma}\sqrt{\frac{2}{\pi}}\exp\left(-\frac{1}{2\sigma^2}\right)~\mbox{as}~x\to\infty,\\
\notag g_2(x)&=1+\frac{\sigma^4}{x^2}\to 1~\mbox{as}~x\to\infty,
\end{align}
and
\begin{align}\label{eq.10}
\tilde{f}_2(x)=2\sigma^2\left(1+\frac{\sigma^2}{x}\right)~\mbox{with}~\tilde{f}_2'(x)\to 0~\mbox{as}~x\to\infty.
\end{align}
\end{lemma}
\noindent \textbf{Proof.}
 Similar to the case of $t\neq 2$, we get
\begin{align} \notag1-F_2(x)&=2\frac{\sigma^2f(x^{\frac{1}{2}})}{x^{\frac{1}{2}}}\left[1+\sigma^2x^{-1}
-\sigma^4x^{-2}+3\sigma^6x^{-3}+O(x^{-4})\right]\\
\label{eq.11}&=2\frac{\sigma^2f(x^{\frac{1}{2}})}{x^{\frac{1}{2}}}(1+\sigma^2x^{-1})\left[1
-\sigma^4x^{-2}(1+\sigma^2x^{-1})^{-1}+3\sigma^6x^{-3}(1+\sigma^2x^{-1})^{-1}+O(x^{-4})\right]\\
\notag&=\frac{2}{\sigma}\sqrt{\frac{2}{\pi}}\exp\left[-\frac{x}{2\sigma^2}+\frac{1}{2}\log x+\log\left(1+\frac{\sigma^2}{x}\right)\right]\left[1-\sigma^4x^{-2}
+4\sigma^6x^{-3}+O(x^{-4})\right]\\
\notag&=\frac{2}{\sigma}\sqrt{\frac{2}{\pi}}\exp\left(-\frac{1}{2\sigma^2}\right)\exp\left(-\int^x_1\frac{g_2(u)}{\tilde{f}_2(u)}du\right)\left[1-\sigma^4x^{-2}
+4\sigma^6x^{-3}+O(x^{-4})\right]
\end{align}
with $g_2(x)=1+\sigma^4t^{-2}$ and $\tilde{f}_2(x)=2\sigma^2(1+\sigma^2t^{-1})$, where the third equality follows from the fact that $(1+x)^a=1+ax+(a(a-1)/2)x^2+O(x^3)$ for all $a\in\mathbb{R}$, as $x\to 0$.

 \qed

Similar to the case of $t\neq 2$, we have the following result:
\begin{proposition}\label{pr.2}
Under the assumptions of Lemma $\ref{le.3}$, we get $F_2(x)\in D(\Lambda)$, where $D(\Lambda)$ is the domain of $\Lambda(x)=\exp\{-\exp(-x)\}.$
\end{proposition}

Now we discuss how to find the constants $c_n,~d_n$. Analogous to the case of $t\neq 2$, we may make choice of $d_n=b^2_n$ and $c_n=\tilde{f}_2(d_n)=2\sigma^2(1+\sigma^2b^{-2}_n)$. Inspired by $c_n$, now change
\begin{align}\label{eq.13}
\notag\bar{d}_n&=b^2_n+2\sigma^4b^{-2}_n,\\
\notag\bar{c}_n&=\tilde{f}_2(\bar{d}_n)\\
\notag&=2\sigma^2[1+\sigma^2b^{-2}_n-2\sigma^6b^{-6}_n+O(b^{-10}_n)]\\
&\sim2\sigma^2(1+\sigma^2b^{-2}_n).
\end{align}

Let
$$T_n(x,t)=F^{n-1}((c_nx+d_n)^{1/t})-(1-F((c_nx+d_n)^{1/t}))^{n-1}.$$
The following lemmas present the expansions of the two terms of densities of $(|M_n|^t-d_n)/c_n$.

\begin{lemma}\label{le.4}
For normalizing constants $c_n$ and $d_n$ determined by (\ref{eq.12}) and $0<t\neq 2$, we have
\begin{align}\label{eq.14}
T_n(x,t)=\Lambda(x)\left\{1-A_1(t,x)e^{-x}b^{-2}_n+\left(\frac{1}{2}A^2_1(t,x)e^{-x}
-A_2(t,x)\right)e^{-x}b^{-4}_n
+O(b^{-6}_n)\right\}
\end{align}
as $n\to\infty$, where
\begin{align}\label{eq.15}
A_1(t,x)=\sigma^2\left(1+x+\frac{(t-2)x^2}{2}\right)
\end{align}
and
\begin{align}\label{eq.16}
A_2(t,x)=\sigma^4\left(\frac{(t-2)^2x^4}{8}+\frac{1}{6}(t-2)(5-2t)x^3-\frac{x^2}{2}-x-1\right).
\end{align}
\end{lemma}
\noindent \textbf{Proof.} Let $\delta_n(x,t)=(c_nx+d_n)^{1/t}.$ One can easily see that $c_nx+d_n>0$ for large $n$ and fixed $x\in\mathbb{R}$. By (\ref{eq.4}), for large $n$, we have $b^2_n\sim2\sigma^2\log n$. Then, by (\ref{eq.12}), we have
\begin{align}\label{eq.17}
\delta^a_n(x,t)=b^a_n\left[1+\frac{a\sigma^2x}{b^2_n}+\frac{a(a-t)\sigma^4x^2}{2b^4_n}
+\frac{a(a-t)(a-2t)\sigma^6x^3}{6b^6_n}+O(b^{-8}_n)\right],
\end{align}
where it follows from the fact that $$(1+x)^a=1+ax+(a(a-1)/2)x^2+(a(a-1)(a-2)/6)x^3+O(x^4),$$ for $a\in\mathbb{R}$, as $x\to 0$. Then, we get
\begin{align}
\notag\frac{\sigma^2 f(\delta_n(x,t))}{\delta_n(x,t)}&\overset{\text{(a)}}{=}\frac{1}{\sigma}\sqrt{\frac{2}{\pi}}b_n\left[1
+\frac{\sigma^2x}{b^2_n}+\frac{(1-t)\sigma^4x^2}{2b^4_n}
+\frac{(1-t)(1-2t)\sigma^6x^3}{6b^6_n}+O(b^{-8}_n)\right]\\
\notag&\times\exp\left\{-\frac{b^2_n}{2\sigma^2}\left[1+\frac{2\sigma^2x}{b^2_n}
+\frac{(2-t)\sigma^4x^2}{b^4_n}
+\frac{(2-t)(2-2t)\sigma^6x^3}{3b^6_n}+O(b^{-8}_n)\right]\right\}\\
\notag&\overset{\text{(b)}}{=}\frac{\sigma^2f(b_n)}{b_n}e^{-x}\left[1+\frac{\sigma^2x}{b^2_n}
+\frac{(1-t)\sigma^4x^2}{2b^4_n}+O(b^{-6}_n)\right]\\
\notag&\times\left[1-\frac{(2-t)\sigma^2x^2}{2b^2_n}-\frac{(2-t)(1-t)\sigma^4x^3}{3b^4_n}
+\frac{(2-t)^2\sigma^4x^4}{8b^4_n}+O(b^{-6}_n)\right]\\
\notag&\overset{\text{(c)}}{=}n^{-1}e^{-x}\bigg\{1+\frac{\sigma^2x}{b^2_n}\left(1+\frac{1}{2}(t-2)x\right)\\
\label{eq.18}&+\frac{\sigma^4x^2}{b^4_n}\left[\frac{1}{8}(t-2)^2x^2+\frac{1}{6}(t-2)(5-2t)x+\frac{1-t}{2}
\right]+O(b^{-6}_n)\bigg\}
\end{align}
where (a) follows from (\ref{eq.17}) with $a=1$ and $2$, (b) is from the fact that $e^x=1+x+x^2/2+O(x^3)$, as $x\to0$ and (c) is due to (\ref{eq.4}). Furthermore, we get
\begin{align}
\notag&1+\sigma^2\delta^{-2}_n(x,t)-\sigma^4\delta^{-4}_n(x,t)+O(\delta^{-6}_n(x,t))\\
\notag&\overset{\text{(a)}}{=}1+\frac{\sigma^2}{b^2_n}\left[1-\frac{2\sigma^2x}{b^2_n}+O(b^{-4}_n)\right]-
\frac{\sigma^4}{b^4_n}\left[1+O(b^{-2}_n)\right]+O(b^{-6}_n)\\
\label{eq.19}&=1+\frac{\sigma^2}{b^2_n}-\frac{\sigma^4}{b^4_n}(1+2x)+O(b^{-6}_n),
\end{align}
where (a) is from (\ref{eq.17}) with $a=-2$ and $-4$. By Lemma $\ref{le.1}$, we get
\begin{align}
\notag1-F(\delta_n(x,t))&=\frac{\sigma^2f(\delta_n(x,t))}{\delta_n(x,t)}\left[1+
\sigma^2\delta^{-2}_n(x,t)-\sigma^4\delta^{-2}_n(x,t)+O(\delta^{-6}_n(x,t))\right]\\
\notag&\overset{\text{(a)}}{=}n^{-1}e^{-x}\bigg\{1+\frac{\sigma^2}{b^2_n}\left[1+x+\frac{1}{2}(t-2)x^2\right]\\
\notag&+\frac{\sigma^4}{b^4_n}\left[\frac{1}{8}(t-2)^2x^4+\frac{1}{6}(t-2)(5-2t)x^3
-\frac{x^2}{2}-x-1\right]+O(b^{-6}_n)\bigg\}\\
\label{eq.20}&=:n^{-1}e^{-x}\left[1+A_1(t,x)b^{-2}_n+A_2(t,x)b^{-4}_n+O(b^{-6}_n)\right],
\end{align}
where (a) is due to (\ref{eq.18}) and (\ref{eq.19}). Accordingly,
\begin{align}
\notag F^{n-1}(\delta_n(x,t))&=\exp\left\{(n-1)\log[1-(1-F(\delta_n(x,t)))]\right\}\\
\notag&\overset{\text{(a)}}{=}\Lambda(x)\exp\left[-A_1(t,x)e^{-x}b^{-2}_n-A_2(t,x)e^{-x}b^{-4}_n+O(b^{-6}_n)\right]\\
\label{eq.21}&\overset{\text{(b)}}{=}\Lambda(x)\left\{1-A_1(t,x)e^{-x}b^{-2}_n+\left(\frac{1}{2}A^2_1(t,x)e^{-x}
-A_2(t,x)\right)e^{-x}b^{-4}_n
+O(b^{-6}_n)\right\},
\end{align}
and
\begin{align}\label{eq.22}
(1-F(\delta_n(x,t)))^{n-1}=\left\{\frac{e^{-x}}{n}\left[1+O(b^{-2}_n)\right]\right\}^{n-1}=o(b^{-\eta}_n),
~\eta\geq6,
\end{align}
where (a) is from the fact that $\log(1-x)=-x+O(x^2)$, as $x\to0$, and (b) follows from that Taylor's expansion of $e^x$. The desired result follows by (\ref{eq.21}) and (\ref{eq.22}).

\qed

\begin{lemma}\label{le.5}
For the normalizing constants $c_n$ and $d_n$ determined by (\ref{eq.12}) and $0<t\neq 2$, we have
\begin{align}\label{eq.23}
\notag n\frac{d}{dx}F((c_nx+d_n)^{1/t})&=e^{-x}\bigg\{1+\frac{\sigma^2x}{b^2_n}[3-t-(2-t)x]\\
&+\frac{\sigma^4x^2}{b^4_n}\left[\frac{1}{2}(3-t)(3-2t)+(t-2)\left(\frac{11}{6}-\frac{5}{6}t\right)x
+\frac{1}{8}(t-2)^2x^2\right]+O(b^{-6}_n)\bigg\},
\end{align}
as $n\to\infty$.
\end{lemma}

\noindent \textbf{Proof.}
It is not hard to check that $$n\frac{d}{dx}F((c_nx+d_n)^{1/t})=\frac{1}{t}nc_n(c_nx+d_n)^{1/t-1}f((c_nx+d_n)^{1/t}).$$
Therefore, we get
\begin{align}
\notag n\frac{d}{dx}F((c_nx+d_n)^{1/t})&\overset{\text{(a)}}{=}\frac{n}{\sigma}b_n\sqrt{\frac{2}{\pi}}\left[1
+\frac{(3-t)\sigma^2x}{b^2_n}+\frac{(3-t)(3-2t)\sigma^4x^2}{2b^4_n}+O(b^{-6}_n)\right]\\
\notag&\times\exp\left\{-\frac{b^2_n}{2\sigma^2}\left[1+\frac{2\sigma^2x}{b^2_n}
+\frac{(2-t)\sigma^4x^2}{b^4_n}+\frac{(2-t)(2-2t)\sigma^6x^3}{3b^6_n}+O(b^{-8}_n)\right]\right\}\\
\notag&\overset{\text{(b)}}{=}nf(b_n)\frac{\sigma^2}{b_n}e^{-x}\left[1
+\frac{(3-t)\sigma^2x}{b^2_n}+\frac{(3-t)(3-2t)\sigma^4x^2}{2b^4_n}+O(b^{-6}_n)\right]\\
\notag&\times\left[1-\frac{(2-t)\sigma^2x^2}{2b^2_n}-\frac{(2-t)(1-t)\sigma^4x^3}{3b^4_n}
+\frac{(2-t)^2\sigma^4x^4}{8b^4_n}+O(b^{-6}_n)\right]\\
\notag&\overset{\text{(c)}}{=}e^{-x}\bigg\{1+\frac{\sigma^2x}{b^2_n}[3-t-(2-t)x]\\
\notag&+\frac{\sigma^4x^2}{b^4_n}\left[\frac{1}{2}(3-t)(3-2t)+(t-2)\left(\frac{11}{6}-\frac{5}{6}t\right)x
+\frac{1}{8}(t-2)^2x^2\right]+O(b^{-6}_n)\bigg\},
\end{align}
where (a) follows from (\ref{eq.17}) with $a=3-t$, (\ref{eq.12}) and (\ref{eq.18}) for the expansion of $f(\delta_n(x,t))$ with $\delta_n(x,t)=(c_nx+d_n)^{1/t}$, (b) is from the fact that $e^x=1+x+x^2/2+O(x^3)$, as $x\to0$ and (c) is due to (\ref{eq.4}). The proof is complete.

 \qed

\begin{lemma}\label{le.6}
For the normalizing constants $c_n$ and $d_n$ determined by (\ref{eq.13}) and $t=2$, we have
\begin{align}\label{eq.24}
T_n(x,t)=\Lambda(x)\left[1-B_1(t,x)e^{-x}b^{-4}_n-B_2(t,x)e^{-x}b^{-6}_n+O(b^{-8}_n)\right],
\end{align}
as $n\to\infty$, where
\begin{align}
\label{eq.25}B_1(t,x)=-\sigma^4\left(x^2+x+\frac{1}{2}\right)
\end{align}
and
\begin{align}
\label{eq.26}B_2(t,x)=\sigma^6\left(\frac{4}{3}x^3+2x^2-2x+\frac{7}{3}\right).
\end{align}
\end{lemma}

\noindent \textbf{Proof.}
The proof of the case of $t=2$ is similar to the case of $0<t\neq2$. Note that $c_n=2\sigma^2(1+\sigma^2b^{-2}_n)$, $d_n=b^2_n+2\sigma^4b^{-2}_n$ for $t=2$. So, we get
$$\delta_n(x,2)=(c_nx+d_n)^{1/2}=b_n[1+2\sigma^2b^{-2}_nx+2\sigma^4(x+1)b^{-4}_n]^{1/2}
=:\beta_n.$$
Then, we have
\begin{align}
\label{eq.27}\beta_n^a=b^a_n\left[1+\frac{a\sigma^2x}{b^2_n}+\frac{a\sigma^4}{b^4_n}\left(1+x
-\frac{2-a}{2}x^2\right)-\frac{a(2-a)\sigma^6x}{b^6_n}\left(1+x
-\frac{4-a}{6}x^2\right)+O(b^{-8}_n)\right].
\end{align}
Further, we get
\begin{align}
\notag\frac{\sigma^2f(\beta_n)}{\beta_n}\overset{\text{(a)}}{=}&\sqrt{\frac{2}{\pi}}\frac{b^2_n}{\sigma^3}
\exp\left(-\frac{b^2_n}{2\sigma^2}\right)\frac{\sigma^2}{b_n}e^{-x}\\
\notag&\times\left[1
+\frac{\sigma^2x}{b^2_n}+\frac{\sigma^4}{b^4_n}
\left(1+x-\frac{1}{2}x^2\right)-\frac{\sigma^6x}{b^6_n}\left(1+x
-\frac{1}{2}x^2\right)+O(b^{-8}_n)\right]\\
\notag&\times\left[1
-\frac{\sigma^2(1+x)}{b^2_n}+\frac{\sigma^4(1+x)^2}{2b^4_n}
-\frac{\sigma^6(1+x)^3}{6b^6_n}+O(b^{-8}_n)\right]\\
\label{eq.28}\overset{\text{(b)}}{=}&n^{-1}e^{-x}\left[1-\frac{\sigma^2}{b^2_n}-
\frac{\sigma^4}{b^4_n}\left(x^2-x-\frac{3}{2}\right)+\frac{\sigma^6}{b^6_n}\left(\frac{4x^3}{3}-x^2-3x-
\frac{7}{6}\right)+O(b^{-8}_n)\right],
\end{align}
where (a) is from (\ref{eq.27}) with $a=1$ and $2$ and $e^x=1+x+x^2/2+O(x^3)$, as $x\to0$, and (b) is due to (\ref{eq.4}). Besides, applying (\ref{eq.27}) with $a=-2$, $-4$ and $-6$, we get
\begin{align}
\notag&1+\sigma^2\beta_n^{-2}-\sigma^4\beta_n^{-4}+3\sigma^6\beta_n^{-6}+O(\beta_n^{-8})\\
\notag=&1+\sigma^2b_n^{-2}\left[1-\frac{2\sigma^2x}{b^2_n}-\frac{2\sigma^4}{b^4_n}\left(1+x
-2x^2\right)+O(b^{-6}_n)\right]\\
\notag&-\sigma^4b_n^{-4}\left[1-\frac{4\sigma^2x}{b^2_n}+O(b^{-4}_n)\right]
+3\sigma^6b_n^{-6}(1+O(b^{-2}_n))+O(b^{-8}_n)\\
\label{eq.29}=&1+\frac{\sigma^2}{b^2_n}-
\frac{\sigma^4}{b^4_n}(2x+1)+\frac{\sigma^6}{b^6_n}(4x^2-2x+1)+O(b^{-8}_n).
\end{align}
Combining with Lemma \ref{le.1}, (\ref{eq.28}) and (\ref{eq.29}), we get
\begin{align}
\notag1-F(\beta_n)&=n^{-1}e^{-x}\left[1-\frac{\sigma^4}{b^4_n}\left(x^2+x+\frac{1}{2}\right)
+\frac{\sigma^6}{b^4_n}\left(\frac{4}{3}x^3+2x^2-2x+\frac{7}{3}\right)+O(b^{-8}_n)\right]\\
\label{eq.30}&=:n^{-1}e^{-x}\left[1+B_1(t,x)b^{-4}_n+B_2(t,x)b^{-6}_n+O(b^{-8}_n)\right].
\end{align}
The remainder proof is the same as the case of $0<t\neq2$. We omit it. The proof is complete.

\qed

\begin{lemma}\label{le.7}
For the normalizing constants $c_n$ and $d_n$ determined by (\ref{eq.13}) and $t=2$, we have
\begin{align}
\notag n\frac{d}{dx}F((c_nx+d_n)^{1/t})&=e^{-x}\bigg\{1-\frac{\sigma^4}{b^4_n}\left(x^2-x-\frac{1}{2}\right)
+\frac{\sigma^6}{b^6_n}\left(\frac{4}{3}x^3-2x^2-2x+\frac{1}{3}\right)+O(b^{-8}_n)\bigg\},
\end{align}
as $n\to\infty$.
\end{lemma}

\noindent \textbf{Proof.}
By (\ref{eq.28}) and after observing that $c_n=2\sigma^2(1+\sigma^2b^{-2}_n)$, we get
\begin{align}\label{eq.31}
\notag n\frac{d}{dx}F(\beta_n)&=e^{-x}\left(1+\frac{\sigma^2}{b^2_n}\right)\left[1-\frac{\sigma^2}{b^2_n}-
\frac{\sigma^4}{b^4_n}\left(x^2-x-\frac{3}{2}\right)+\frac{\sigma^6}{b^6_n}\left(\frac{4x^3}{3}-x^2-3x-
\frac{7}{6}\right)+O(b^{-8}_n)\right]\\
&=e^{-x}\bigg\{1-\frac{\sigma^4}{b^4_n}\left(x^2-x-\frac{1}{2}\right)
+\frac{\sigma^6}{b^6_n}\left(\frac{4}{3}x^3-2x^2-2x+\frac{1}{3}\right)+O(b^{-8}_n)\bigg\}.
\end{align}
The proof is complete.

 \qed

As we mentioned in the introduction, Liu and Liu \cite{Liu and Liu} obtained the pointwise convergence rate of distribution of partial maximum to its limiting distribution. Their main results are stated as follows.

\begin{theorem}\label{th.3}
Suppose that $\{X_n,n\geq1\}$ is a sequence of i.i.d. random variables with cdf MD. Then,
\begin{align}\label{eq.47}
F^n(\hat{a}_nx+\hat{b}_n)-\Lambda(x)\sim\Lambda(x)e^{-x}\frac{(\log(2\log n))^2}{16\log n},
\end{align}
for large $n$, where
\begin{align}\label{eq.48}
\hat{a}_n=\frac{\sigma}{(2\log n)^{1/2}}~\mbox{and}~\hat{b}_n=(2\sigma^2\log n)^{1/2}+\frac{\sigma\log(2\log n)+\sigma\log \frac{2}{\pi}}{2(2\log n)^{1/2}}.
\end{align}
\end{theorem}

\section{Main result}
\label{sec3}

In this section, we establish the higher-order expansions of the cdf and the pdf of
powered maximum from MD sample.
\begin{theorem}\label{th.1}
~\\
(i) For $0<t\neq2$ and the normalizing constants $c_n$ and $d_n$ given by (\ref{eq.12}), we have
\begin{align}\label{eq.32}
\mathbb{P}(|M_n|^t\leq c_nx+d_n)=\Lambda(x)\left\{1-e^{-x}A_1(t,x)b^{-2}_n+e^{-x}\left[\frac{1}{2}e^{-x}A^2_1(t,x)
-A_2(t,x)\right]b^{-4}_n+O(b^{-6}_n)\right\},
\end{align}
where
\begin{align}\label{eq.33}
A_1(t,x)=\sigma^2\left[1+x+\frac{1}{2}(t-2)x^2\right]
\end{align}
and
\begin{align}\label{eq.34}
A_2(t,x)=\sigma^4\left[\frac{1}{8}(t-2)^2x^4+\frac{1}{6}(t-2)(5-2t)x^3
-\frac{x^2}{2}-x-1\right].
\end{align}
~\\
(ii) For $t=2$ and the normalizing constants $c_n$ and $d_n$ given by (\ref{eq.13}), we have
\begin{align}\label{eq.35}
\mathbb{P}(|M_n|^t\leq c_nx+d_n)=\Lambda(x)\left[1-e^{-x}B_1(t,x)b^{-4}_n-e^{-x}B_2(t,x)b^{-6}_n+O(b^{-8}_n)\right],
\end{align}
where
\begin{align}\label{eq.36}
B_1(t,x)=-\sigma^4\left(x^2+x+\frac{1}{2}\right)
\end{align}
and
\begin{align}\label{eq.37}
B_2(t,x)=\sigma^6\left(\frac{4}{3}x^3+2x^2-2x+\frac{7}{3}\right).
\end{align}
\end{theorem}

\begin{remark}\label{re.1}
From Theorem \ref{th.1}, one can easily see that the convergence rates of powered maximum of cdf for MD are proportional to $1/\log n$ and $1/(\log n)^2$ for power index $0<t\neq2$ and $t=2$, respectively, since $1/b^2_n\sim2\sigma^2\log n$ by (\ref{eq.4}).
\end{remark}

\begin{remark}\label{re.5}
From Theorems \ref{th.3} and \ref{th.1} (ii), we can observe that the convergence speed of powered extreme of cdf for MD is better than that of extreme of cdf.
\end{remark}

In the following we provide the higher-order expansions of the pdf of powered maximum.
\begin{theorem}\label{th.2}
~\\
(i) For $0<t\neq2$ and the normalizing constants $c_n$ and $d_n$ given by (\ref{eq.12}), we have
\begin{align}
\label{eq.41}\frac{d}{dx}\mathbb{P}(|M_n|^t\leq c_nx+d_n)=\Lambda'(x)\left[1+P_1(t,x)b^{-2}_n+P_2(t,x)b^{-4}_n+O(b^{-6}_n)\right],
\end{align}
where
\begin{align}
\notag P_1(t,x)=\sigma^2\left\{-\left[\frac{(t-2)x^2}{2}+x+1\right]e^{-x}+(t-2)x^2-(t-3)x\right\}
\end{align}
and
\begin{align}
\notag P_2(t,x)=&\sigma^4\bigg\{\frac{1}{2}\left[\frac{(t-2)x^2}{2}+x+1\right]^2e^{-2x}\\
\notag&-\left[\frac{5(t-2)x^4}{8}-(t-2)\left(\frac{5}{6}t-\frac{10}{3}\right)x^3
+\left(2t+\frac{1}{2}\right)x^2-1\right]e^{-x}\\
\notag&+\frac{(t-2)^2x^3}{8}-(t-2)\left(\frac{5}{6}t-\frac{11}{6}\right)x^2
+\frac{(t-3)(2t-3)}{2}x\bigg\}.
\end{align}
~\\
(ii) For $t=2$ and the normalizing constants $c_n$ and $d_n$ given by (\ref{eq.13}), we have
\begin{align}
\label{eq.42}\frac{d}{dx}\mathbb{P}(|M_n|^t\leq c_nx+d_n)=\Lambda'(x)\left[1+Q_1(t,x)b^{-4}_n+Q_2(t,x)b^{-6}_n+O(b^{-8}_n)\right],
\end{align}
where
\begin{align}
\notag Q_1(t,x)=\sigma^4\left[\left(x^2+x+\frac{1}{2}\right)e^{-x}-x^2+x+\frac{1}{2}\right]
\end{align}
and
\begin{align}
\notag Q_2(t,x)=-\sigma^6\left[\left(\frac{4}{3}x^3+2x^2-2x+\frac{7}{3}\right)e^{-x}-
\frac{4}{3}x^3+2x^2+2x-\frac{1}{3}\right].
\end{align}
\end{theorem}

\begin{remark}\label{re.2}
From Theorem \ref{th.2}, it is not difficult to observe that the convergence speeds of powered extreme of pdf for MD are the same order of $1/\log n$ and $1/(\log n)^2$ for power index $0<t\neq2$ and $t=2$, respectively, because of $1/b^2_n\sim2\sigma^2\log n$ by (\ref{eq.4}).
\end{remark}

\begin{remark}\label{re.3}
For $t=2$, the normalizing constants $c_n$ and $d_n$ are not given by (\ref{eq.13}), but we choose them as follows:
\begin{align}\label{eq.44}
c_n=2\sigma^2(1-\sigma^2b^{-2}_n)~\mbox{and}~d_n=b^2_n-2\sigma^4b^{-2}_n,
\end{align}
then we derive
\begin{align}
\notag&\mathbb{P}(|M_n|^t\leq c_nx+d_n)\\
\notag=&\Lambda(x)\bigg\{1-\frac{2e^{-x}\sigma^2}{b^{2}_n}(x+1)+\frac{e^{-x}\sigma^4}{b^{4}_n}\left[2e^{-x}(x+1)^2
-x^2-x-\frac{3}{2}\right]b^{-4}_n\\
\notag&-\frac{e^{-x}\sigma^6}{b^{6}_n}\bigg[\frac{4}{3}e^{-2x}(x+1)^3
-2e^{-x}(x+1)\left(x^2+x+\frac{3}{2}\right)\\
\label{eq.45}&+\frac{2}{3}x^3+2x^2+3x+\frac{14}{3}\bigg]+O(b^{-8}_n)\bigg\}
\end{align}
and
\begin{align}
\notag&\frac{d}{dx}\mathbb{P}(|M_n|^t\leq c_nx+d_n)\\
\notag&=\Lambda'(x)\bigg\{1-\frac{2\sigma^2}{b^{2}_n}[e^{-x}(x+1)-x]
+\frac{\sigma^4}{b^{4}_n}\bigg[2e^{-2x}(x+1)^2-\left(5x^2+5x+\frac{3}{2}\right)e^{-x}\\
\label{eq.46}&+x^2-x+\frac{1}{2}\bigg]+\frac{\sigma^6}{b^{6}_n}\bigg[4x(x+1)^2e^{-2x}
-(4x^3+2x^2+2x+1)e^{-x}+\frac{2}{3}x^3-x-\frac{7}{6}\bigg]+O(b^{-8}_n)\bigg\}.
\end{align}
Obviously, the convergence rates of the cdf and the pdf of powered extreme given by (\ref{eq.35}) and (\ref{eq.42}), which are proportional to $1/(\log n)^2$, are faster than that given by (\ref{eq.45}) and (\ref{eq.46}). Consequently, the normalizing constants $c_n$ and $d_n$ determined by (\ref{eq.13}) are optimal.
\end{remark}

\section{Numerical analysis}
In this section, we conduct numerical studies to illustrate the accurateness of higher-order expansions for the cdf and the pdf of $|M_n|^t$. Let $T^{(i)}(x)$ and $S^{(i)}(x)$, $i=1,2,3,$ respectively represent the first-order, the second-order and the third-order approximations of the cdf and the pdf of $|M_n|^t$. Since the analysis of the case of $t\neq2$ is similar to that of $t=2$, we only consider the situation of $t=2$. By Theorems \ref{th.1} and \ref{th.2}, we obtain
\begin{align}
\notag&T^{(1)}(x)=\Lambda(x),\\
\notag&T^{(2)}(x)=\Lambda(x)\left[1-e^{-x}B_1(t,x)b^{-4}_n\right],\\
\notag&T^{(3)}(x)=\Lambda(x)\left[1-e^{-x}B_1(t,x)b^{-4}_n-e^{-x}B_2(t,x)b^{-6}_n\right],
\end{align}
and
\begin{align}
\notag&S^{(1)}(x)=\Lambda(x)\exp(-x),\\
\notag&S^{(2)}(x)=\Lambda(x)\exp(-x)\left[1+Q_1(t,x)b^{-4}_n\right],\\
\notag&S^{(3)}(x)=\Lambda(x)\exp(-x)\left[1+Q_1(t,x)b^{-4}_n+Q_2(t,x)b^{-6}_n\right].
\end{align}
Easily observe that the second-order approximation and the third-order relate to the sample size $n$.

In order to compare the precision of true values with its approximations, let
\begin{align}
\notag&E^{(i)}(x)=\left|F^n(\sqrt{c_nx+d_n})-T^{(i)}(x)\right|
\end{align}
and
\begin{align}
\notag&G^{(i)}(x)=\left|\frac{nc_n}{2\sqrt{c_nx+d_n}}F^{n-1}(\sqrt{c_nx+d_n})f(\sqrt{c_nx+d_n})-S^{(i)}(x)\right|
\end{align}
respectively stand for the absolute errors of the cdf and the pdf, where $i=1,2,3.$ We utilize MATLAB to compute the approximations and the true values of the cdf and the pdf of $M^2_n$.

First, we estimate the absolute errors of the cdf of $M^2_n$ at $x=0.7$, where the sample size $n$ varies from $25$ to $1000$ with step size $25$. For given $x=0.7$, numerical analysis results of $E^{(i)}(x)$ are recorded in Table \ref{tab.1}. The table demonstrates that the precision of all three kinds of approximations of the cdf can be refined as the sample size $n$ increases.

\begin{figure}[h]
\begin{center}
\subfigure[$\sigma=2$]{\includegraphics[width=0.40\textwidth]{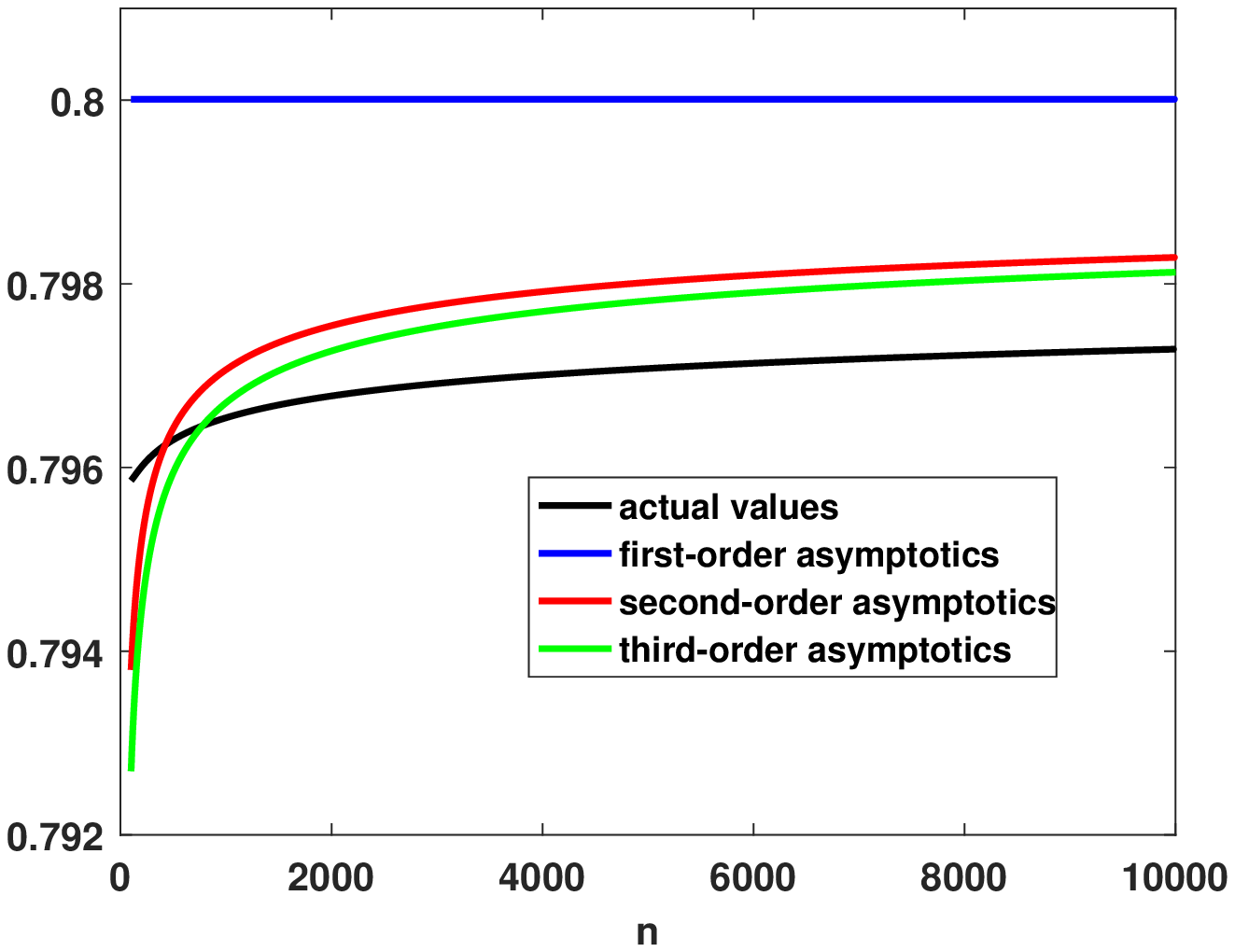}}
\subfigure[$\sigma=0.5$]{\includegraphics[width=0.40\textwidth]{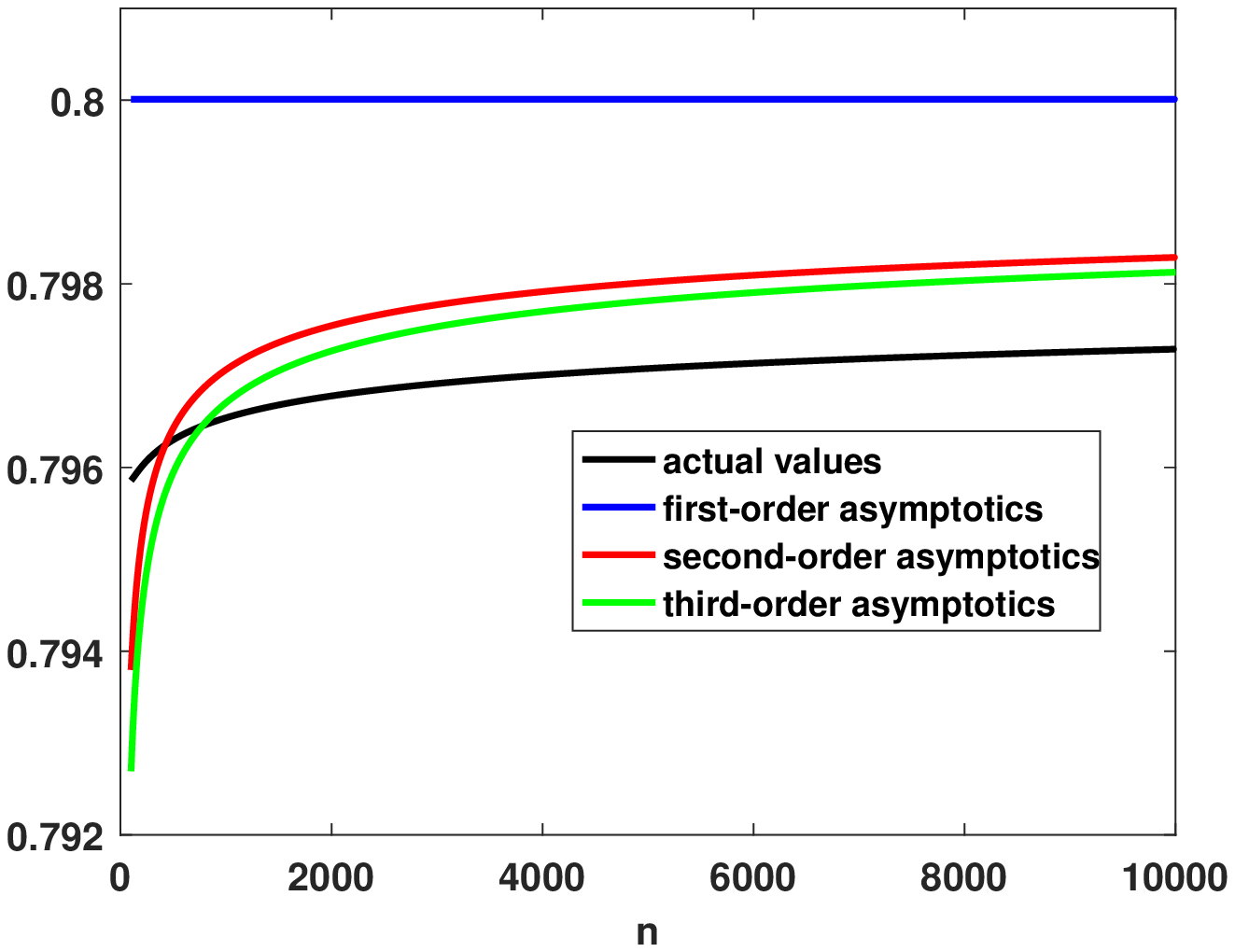}}
\caption{Actual values and its asymptotics of the cdf of $M^2_n$ with $x=1.5$. The actual values drawn in black, the first-order approximations drawn in blue, the second-order approximations drawn in red and the third-order approximation drawn in green.}\label{fig.2}
\end{center}
\vspace*{-14pt}
\end{figure}

To order to indicate the precision of all approximations more intuitive with the change of the sample size $n$, the actual values and its approximation of the cdf of $M^2_n$ are plotted versus the values of $n$ with $x=1.5$. Figure \ref{fig.2} evidences that the larger $n$, the better all asymptotics.

Secondly, we estimate the absolute errors of the pdf of $M^2_n$ at $x=0.7$, where the value of the sample size $n$ ranges from $375$ to $15000$ with step length $375$. Table \ref{tab.2} lists the numerical analysis results of $G^{(i)}(x)$, where $i=1,2,3.$ Table \ref{tab.2} reveals that the precision of all three kinds of approximations of the pdf enhances as the sample size $n$ grows.

\begin{figure}[h]
\begin{center}
\subfigure[$\sigma=2$]{\includegraphics[width=0.40\textwidth]{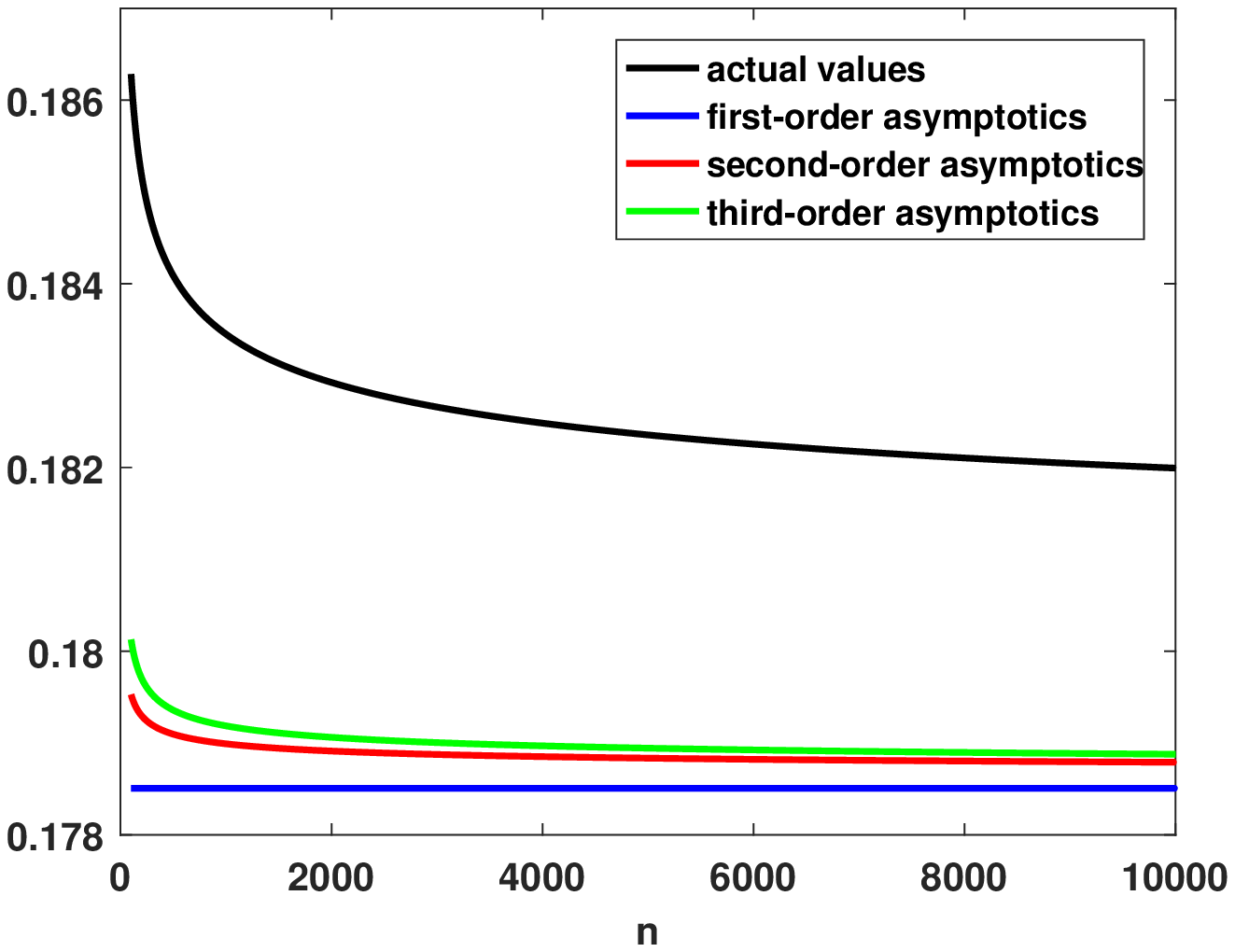}}
\subfigure[$\sigma=0.5$]{\includegraphics[width=0.40\textwidth]{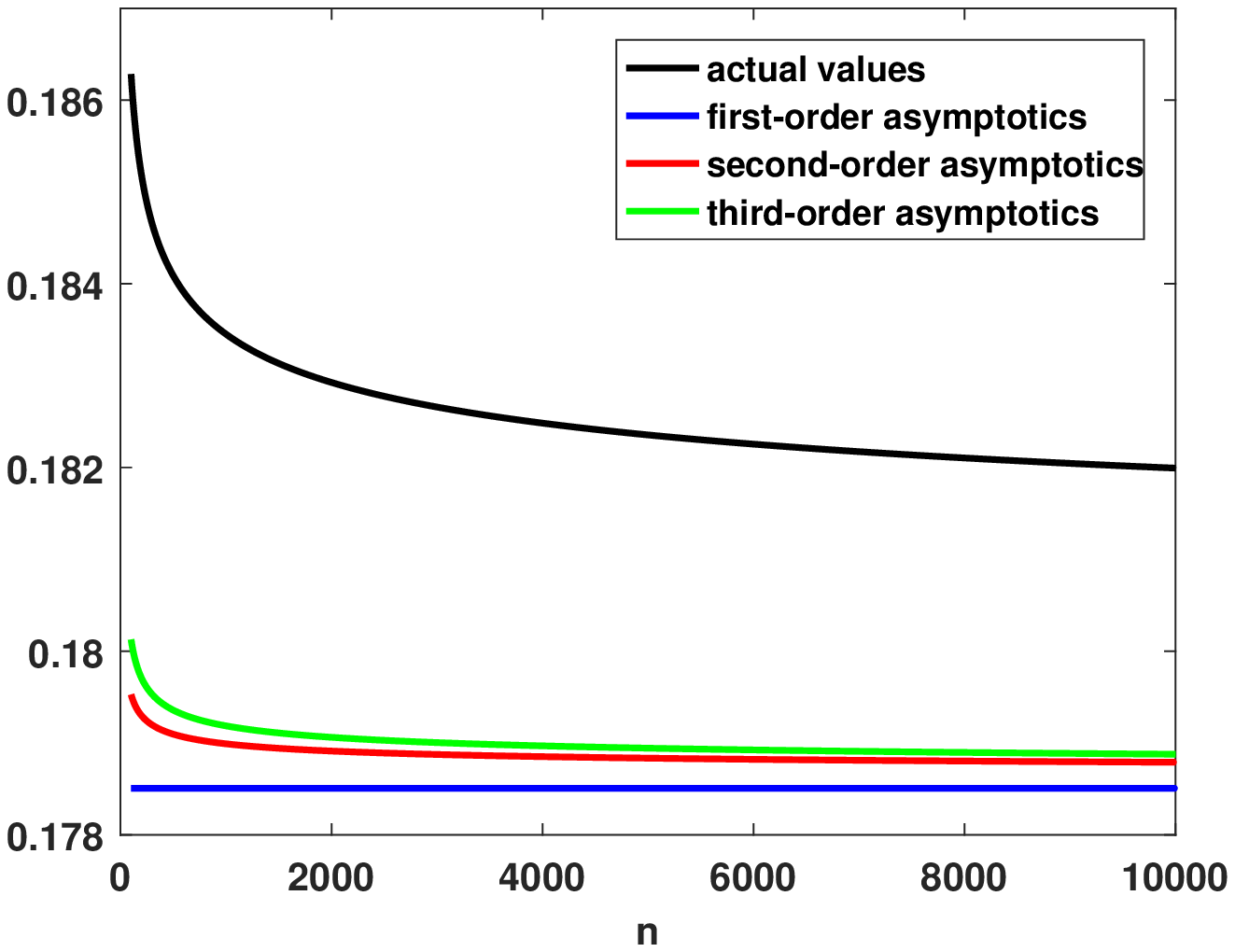}}
\caption{Actual values and its asymptotics of the pdf of $M^2_n$ with $x=1.5$. The actual values drawn in black, the first-order approximations drawn in blue, the second-order approximations drawn in red and the third-order approximation drawn in green.}\label{fig.3}
\end{center}
\vspace*{-14pt}
\end{figure}

To clear the precision of all approximations more intuitive with $n$, the actual and its approximations of the pdf of $M^2_n$ are plotted versus the values of $n$ with $x=1.5$. Figure \ref{fig.3} indicates that as the sample size $n$ becomes larger, all approximations become better.

\begin{table}[htbp]
  \centering
  \large
    \begin{tabular}{rrrr}
    \hline
         $n$ & $E^{(1)}(x)$ & $E^{(2)}(x)$ & $E^{(3)}(x)$ \\
    \hline
         25 & 0.0169056391 & 0.00877452615 & 0.00733539417 \\
         50 & 0.0143357459 & 0.00869068028 & 0.00785819009 \\
         75 & 0.0131346277 & 0.00843346219 & 0.00780078242 \\
         100 & 0.0123911158 & 0.00821865886 & 0.00768964941 \\
         125 & 0.0118668421 & 0.00804347997 & 0.00757945239 \\
         150 & 0.0114683134 & 0.0078976489 & 0.00747885611 \\
         175 & 0.0111502039 & 0.00777354585 & 0.00738841683 \\
         200 & 0.0108874336 & 0.00766594114 & 0.00730705118 \\
         225 & 0.0106648041 & 0.00757120115 & 0.00723346892 \\
         250 & 0.0104724714 & 0.00748673237 & 0.00716650872 \\
         275 & 0.0103037264 & 0.00741063197 & 0.00710519645 \\
         300 & 0.0101538089 & 0.00734146835 & 0.00704873143 \\
         325 & 0.0100192298 & 0.00727814037 & 0.0069964575 \\
         350 & 0.00989736162 & 0.00721978455 & 0.00694783483 \\
         375 & 0.00978618048 & 0.00716571218 & 0.00690241634 \\
         400 & 0.00968409693 &  0.0071153658 & 0.00685982891 \\
         425 & 0.00958984178 & 0.00706828819 & 0.00681975869 \\
         450 & 0.00950238672 & 0.00702409995 & 0.00678193964 \\
         475 & 0.00942088789 & 0.00698248302 & 0.00674614458 \\
         500 & 0.00934464492 & 0.00694316822 & 0.00671217821 \\
         525 & 0.00927307061 & 0.00690592583 & 0.00667987148 \\
         550 & 0.00920566819 & 0.00687055834 & 0.00664907727 \\
         575 & 0.00914201391 & 0.00683689469 & 0.00661966674 \\
         600 & 0.00908174361 & 0.00680478587 & 0.00659152652 \\
         625 & 0.00902454227 & 0.00677410125 & 0.00656455636 \\
         650 & 0.00897013561 & 0.00674472575 & 0.00653866719 \\
         675 & 0.00891828351 & 0.00671655745 & 0.00651377957 \\
         700 & 0.00886877463 & 0.00668950571 & 0.00648982235 \\
         725 & 0.00882142208 & 0.00666348955 & 0.00646673154 \\
         750 & 0.00877605982 & 0.00663843637 & 0.00644444946 \\
         775 & 0.00873253972 & 0.00661428081 & 0.00642292387 \\
         800 & 0.00869072914 & 0.00659096382 & 0.00640210739 \\
         825 & 0.00865050883 & 0.00656843192 & 0.00638195685 \\
         850 & 0.00861177127 & 0.00654663651 & 0.00636243286 \\
         875 & 0.00857441915 & 0.00652553326 & 0.00634349935 \\
         900 & 0.00853836418 & 0.00650508169 & 0.00632512325 \\
         925 & 0.008503526 & 0.0064852447 & 0.00630727414 \\
         950 & 0.00846983127 & 0.00646598823 & 0.00628992398 \\
         975 & 0.0084372129 & 0.00644728091 & 0.00627304687 \\
         1000 & 0.00840560939 & 0.00642909381 & 0.00625661887 \\
    \hline
    \end{tabular}%
  \label{tab.1}%
  \caption{Absolute errors between actual values and their asymptotics of the cdf at $x=0.7$ with $\sigma=2$}
\end{table}%

\begin{table}[htbp]
  \centering
  \large
    \begin{tabular}{rrrr}
    \hline
      $n$ & $G^{(1)}(x)$ & $G^{(2)}(x)$ & $G^{(3)}(x)$ \\
    \hline
      375 & 0.00825613746 & 0.00585394461 & 0.00554667797 \\
      750 & 0.00710011928 & 0.00514055207 & 0.00491416905 \\
      1125 & 0.0065538405 & 0.00479753582 & 0.00460544643 \\
      1500 & 0.00621014157 & 0.00457905959 & 0.00440714319 \\
      1875 & 0.00596472382 & 0.00442157961 & 0.00426337709 \\
      2250 & 0.00577637198 & 0.00429979382 & 0.0041517166 \\
      2625 & 0.00562489953 & 0.00420122856 & 0.00406103823 \\
      3000 & 0.00549902795 & 0.00411887475 & 0.0039850629 \\
      3375 & 0.00539186326 & 0.00404842643 & 0.00391991863 \\
      3750 & 0.00529890627 & 0.00398706069 & 0.00386305898 \\
      4125 & 0.00521707021 & 0.0039328329 & 0.00381272502 \\
      4500 & 0.0051441522 & 0.00388435007 & 0.00376765376 \\
      4875 & 0.00507852912 & 0.00384058231 & 0.0037269095 \\
      5250 & 0.00501897306 & 0.00380074803 & 0.00368978076 \\
      5625 & 0.00496453429 & 0.00376424093 & 0.0036557147 \\
      6000 & 0.00491446415 & 0.00373058177 & 0.00362427365 \\
      6375 & 0.00486816283 & 0.00369938562 & 0.00359510557 \\
      6750 & 0.00482514285 & 0.00367033886 & 0.00356792328 \\
      7125 & 0.00478500308 & 0.00364318286 & 0.00354248964 \\
      7500 & 0.00474740971 & 0.00361770192 & 0.00351860668 \\
      7875 & 0.00471208225 & 0.00359371443 & 0.00349610752 \\
      8250 & 0.00467878285 & 0.00357106611 & 0.00347485024 \\
      8625 & 0.00464730824 & 0.00354962483 & 0.00345471324 \\
      9000 & 0.0046174834 & 0.00352927666 & 0.00343559152 \\
      9375 & 0.00458915666 & 0.00350992269 & 0.00341739387 \\
      9750 & 0.00456219575 & 0.00349147652 & 0.00340004057 \\
      10125 & 0.00453648471 & 0.0034738623 & 0.00338346154 \\
      10500 & 0.00451192133 & 0.00345701305 & 0.0033675949 \\
      10875 & 0.00448841509 & 0.0034408694 & 0.00335238574 \\
      11250 & 0.00446588547 & 0.00342537844 & 0.00333778509 \\
      11625 & 0.00444426057 & 0.00341049288 & 0.00332374917 \\
      12000 & 0.00442347592 & 0.00339617026 & 0.00331023861 \\
      12375 & 0.00440347351 & 0.00338237231 & 0.00329721797 \\
      12750 & 0.00438420096 & 0.00336906445 & 0.00328465516 \\
      13125 & 0.00436561084 & 0.00335621534 & 0.00327252109 \\
      13500 & 0.00434766007 & 0.00334379647 & 0.00326078929 \\
      13875 & 0.00433030941 & 0.00333178183 & 0.0032494356 \\
      14250 & 0.00431352301 & 0.00332014766 & 0.00323843795 \\
      14625 & 0.00429726807 & 0.00330887219 & 0.00322777608 \\
      15000 & 0.00428151447 & 0.00329793539 & 0.00321743138 \\
    \hline
    \end{tabular}%
  \label{tab.2}%
  \caption{Absolute errors between actual values and their asymptotics of the pdf at $x=0.7$ with $\sigma=2$}
\end{table}%

\section{Proof of main result}
\noindent \textbf{Proof of Theorem \ref{th.1}.}
By some fundamental calculations, we get
\begin{align}
\label{eq.38}\mathbb{P}(|M_n|^t\leq c_nx+d_n)=F^n((c_nx+d_n)^{1/t})-(1-F((c_nx+d_n)^{1/t}))^n.
\end{align}
First, we consider the case of $0<t\neq2$. By (\ref{eq.20}) and similar discussions as for (\ref{eq.21}) and (\ref{eq.22}), we get
\begin{align}
\label{eq.39} F^n(\delta_n(x,t))=\Lambda(x)\left\{1-A_1(t,x)e^{-x}b^{-2}_n+\left(\frac{1}{2}A^2_1(t,x)e^{-x}
-A_2(t,x)\right)e^{-x}b^{-4}_n
+O(b^{-6}_n)\right\},
\end{align}
where $A_1(t,x)$ and $A_2(t,x)$ are determined by (\ref{eq.15}) and (\ref{eq.16}), and
\begin{align}
\label{eq.40} (1-F(\delta_n(x,t)))^n=\left\{\frac{e^{-x}}{n}\left[1+O(b^{-2}_n)\right]\right\}^{n}=o(b^{-\eta}_n),
~\eta\geq6.
\end{align}
A combination of (\ref{eq.39}) and (\ref{eq.40}) implies that (\ref{eq.32}) holds.

 For the case of $t=2$, by similar arguments as for $0<t\neq2$, the desired result follows. The proof is complete.

\qed

\noindent \textbf{Proof of Theorem \ref{th.2}.} One can easily check that
\begin{align}
\label{eq.43}\frac{d}{dx}\mathbb{P}(|M_n|^t\leq c_nx+d_n)=&n\left(\frac{d}{dx}F((c_nx+d_n)^{1/t})\right)\\
\notag&\times\left\{F^{n-1}((c_nx+d_n)^{1/t})
+[1-F((c_nx+d_n)^{1/t})]^{n-1}\right\}.
\end{align}
For $0<t\neq2$, combining with Lemmas \ref{le.4} and \ref{le.5}, we get
\begin{align}
\notag\frac{1}{\Lambda'(x)}&\frac{d}{dx}\mathbb{P}(|M_n|^t\leq c_nx+d_n)-1=\bigg\{1+\frac{\sigma^2x}{b^2_n}[3-t-(2-t)x]\\
\notag&+\frac{\sigma^4x^2}{b^4_n}\left[\frac{1}{2}(3-t)(3-2t)+(t-2)\left(\frac{11}{6}-\frac{5}{6}t\right)x
+\frac{1}{8}(t-2)^2x^2\right]+O(b^{-6}_n)\bigg\}\\
\notag&\times\bigg\{1-\sigma^2\left[1+x+\frac{1}{2}(t-2)x^2\right]e^{-x}b^{-2}_n
+\bigg(\frac{1}{2}\left[1+x+\frac{1}{2}(t-2)x^2\right]^2e^{-x}\\
\notag&-\left[\frac{1}{8}(t-2)^2x^4+\frac{1}{6}(t-2)(5-2t)x^3
-\frac{x^2}{2}-x-1\right]\bigg)\sigma^4e^{-x}b^{-4}_n
+O(b^{-6}_n)\bigg\}-1\\
\notag=&\frac{\sigma^2}{b^2_n}\left\{-\left[\frac{(t-2)x^2}{2}+x+1\right]e^{-x}+(t-2)x^2-(t-3)x\right\}\\
\notag&+\frac{\sigma^4}{b^4_n}\bigg\{\frac{1}{2}\left[\frac{(t-2)x^2}{2}+x+1\right]^2e^{-2x}\\
\notag&-\left[\frac{5(t-2)x^4}{8}-(t-2)\left(\frac{5}{6}t-\frac{10}{3}\right)x^3
+\left(2t+\frac{1}{2}\right)x^2-1\right]e^{-x}\\
\notag&+\frac{(t-2)^2x^3}{8}-(t-2)\left(\frac{5}{6}t-\frac{11}{6}\right)x^2
+\frac{(t-3)(2t-3)}{2}x\bigg\}+O(b^{-6}_n)\\
\notag=&P_1(t,x)b^{-2}_n+P_2(t,x)b^{-4}_n+O(b^{-6}_n),
\end{align}
which deduces (\ref{eq.41}).

The following is for the case of $t=2$. By (\ref{eq.43}) and Lemmas \ref{le.6} and \ref{le.7}, we gain
\begin{align}
\notag\frac{1}{\Lambda'(x)}&\frac{d}{dx}\mathbb{P}(|M_n|^t\leq c_nx+d_n)-1=\bigg\{1-\frac{\sigma^4}{b^4_n}\left(x^2-x-\frac{1}{2}\right)\\
\notag&+\frac{\sigma^6}{b^6_n}\left(\frac{4}{3}x^3-2x^2-2x+\frac{1}{3}\right)+O(b^{-8}_n)\bigg\}\\
\notag&\times\left[1+\frac{\sigma^4e^{-x}}{b^{4}_n}\left(x^2+x+\frac{1}{2}\right)
-\frac{\sigma^6e^{-x}}{b^{6}_n}\left(\frac{4}{3}x^3+2x^2-2x+\frac{7}{3}\right)+O(b^{-8}_n)\right]-1\\
\notag=&\frac{\sigma^4}{b^{4}_n}\left[\left(x^2+x+\frac{1}{2}\right)e^{-x}-x^2+x+\frac{1}{2}\right]\\
\notag&-\frac{\sigma^6}{b^{6}_n}\left[\left(\frac{4}{3}x^3+2x^2-2x+\frac{7}{3}\right)e^{-x}-
\frac{4}{3}x^3+2x^2+2x-\frac{1}{3}\right]+O(b^{-8}_n)\\
\notag=&Q_1(t,x)b^{-4}_n+Q_2(t,x)b^{-6}_n+O(b^{-8}_n),
\end{align}
which proves (\ref{eq.42}). The proof of Theorem \ref{th.2} is finished.

\qed

\noindent {\bf Acknowledgments}


\noindent {\bf Funding}

This work was supported by Natural Science Foundation of China [grant number 61673015], [grant number 61273020] and Fundamental Research Funds for the Central
Universities [grant number XDJK2015A007], Youth Science and technology talent development project (No.Qian jiao he KY zi [2018]313), Science and technology Foundation of Guizhou province [grant number Qian ke he Ji Chu [2016]1161], Guizhou province natural science foundation in China [grant number Qian Jiao He KY [2016]255].


\begin{thebibliography}{100} \small

\bibitem{Smith}
Smith, LR: Uniform rates of convergence in extreme-value theory. Adv. Appl. Probab. 14, 600-622 (1982)

\bibitem{Leadbetter et al.}
Leadbetter, MR, Lindgren, G, Rootz$\acute{e}$n, H: Extremes and Related Properties of Random Sequences and Processes.
Springer, New York (1983)

\bibitem{Galambos}
Galambos, J. (1987). {\it The asympotic theory of extreme order statistics.} (Second Edition) New York, Wiley.

\bibitem{de Haan and Resnick}
de Haan, L, Resnick, SI: Second-order regular variation and rates of convergence in extreme-value theory. Ann.
Probab. 1, 97-124 (1996)

\bibitem{Hall 1979}
Hall, P., 1979. On the rate of convergence of normal extremes. J. Appl. Probab. 16, 433-439.

\bibitem{Hall 1980}
Hall, P., 1980. Estimating probabilities for normal extremes. Adv. Appl. Probab. 12, 491-500.

\bibitem{Nair}
Nair, K. A. (1981). Asymptotic distribution and moments of normal extremes. Annals of Probability, 9, 150-153.

\bibitem{Liao and Peng}
Liao, X., \& Peng, Z. (2012). Convergence rates of limit distribution of maxima of lognormal samples. Journal of Mathematical Analysis and Applications, 395,
643-653.

\bibitem{Lin et al 2011}
Lin, F., Zhang, X., Peng, Z., \& Jiang, Y. (2011). On the rate of convergence of stsd extremes. Communications in Statistics - Theory and Methods, 40(10), 1795-1806.

\bibitem{Lin et al 2016}
Lin, F., Peng, Z., \& Yu, K. (2016). Convergence rate of extremes for the generalized short-tailed symmetric distribution. Bulletin of the Korean Mathematical Society, 53(5), 1549-1566.

\bibitem{Du and Chen}
Du, L., \& Chen, S. (2016). Asymptotic properties for distributions and densities of extremes from generalized gamma distribution. Journal of the Korean Statistical Society, 45(2), 188-198.

\bibitem{Chen and Du}
Chen, S., \& Du, L. (2017). Asymptotic expansions of density of normalized extremes from logarithmic general error distribution. Communications in Statistics, 46(7), 3459-3478.

\bibitem{Huang et al 2017}
Huang, J., Wang, J., \& Luo, G. (2017a). On the rate of convergence of maxima for the generalized maxwell distribution. Statistics A Journal of Theoretical \& Applied Statistics, 1-13.

\bibitem{Zhou and Ling}
Zhou, W., Ling, C. (2016). Higher-order expansions of powered extremes of normal samples. Statistics and Probability Letters, 111, 12-17.

\bibitem{Liao et al}
Liao, X, Peng, Z, Nadarajah, S: Asymptotic expansions of the moments of skew-normal extremes. Stat. Probab. Lett.
83, 1321-1329 (2013)

\bibitem{Jia et al}
Jia, P, Liao, X, Peng, Z: Asymptotic expansions of the moments of extremes from general error distribution. J. Math.
Anal. Appl. 422, 1131-1145 (2015)

\bibitem{Mandl}
Mandl, F. (2008). {\it Statistical Physics (2nd Edition)}. New Jersey: John Wiley \& Sons.

\bibitem{Shim and Gatignol}
Shim, J. W., \& Gatignol, R. (2013). How to obtain higher-order multivariate hermite expansion of maxwell¨Cboltzmann distribution by using taylor expansion? Zeitschrift F¨¹r Angewandte Mathematik Und Physik, 64(3), 473-482.

\bibitem{Tomer and Panwar}
Tomer, S.K. and Panwar, M.S. (2015). Estimation procedures for Maxwell distribution under type-I progressive hybrid censoring scheme. Journal of Statistical Computation and Simulation, 85(2), 339-356.

\bibitem{Shim}
Shim, J. W. (2017). Parametric lattice boltzmann method. Journal of Computational Physics, 338.

\bibitem{Liu and Liu}
Liu, C., Liu, B. (2013). Convergence rate of extremes from Maxwell sample. J. Inequal. Appl. 2013: 477.
Available at: http://www.journalofinequalitiesandapplications.com/content/2013/1/477.

\bibitem{Huang and Chen}
Huang, J., Chen, S. Tail behavior of the generalized Maxwell distribution. Communications in Statistics-Theory and Methods. 2016, 45(14): 4230-4236.

\bibitem{Dar et al}
Dar, A.A., A. Ahmed and J.A. Reshi. Bayesian analysis of Maxwell-Boltzmann distribution under different loss functions and prior distributions. Pak. J. Statist.
2017, 33(6), 419-440

\bibitem{Huang et al 2017b}
Huang, J., Wang, J., Luo, G., \& He, J. Tail properties and approximate distribution and expansion for extreme of LGMD. Journal of Inequalities \& Applications, 2017b, 2017(1):1-16.

\bibitem{Huang et al 2018}
Huang, J., Wang, J., Luo, G. Pu, H. Higher order expansion for moments of extreme for generalized Maxwell distribution, Communications in Statistics-Theory and Methods, 2018, 47(14): 3441-3452.

\bibitem{Huang and Wang 2018a}
Huang J.W., Wang J.J. On asymptotic of extremes from generalized Maxwell distribution. Bulletin of the Korean Mathematical Society, 2018a, 55(3): 679-698.

\bibitem{Huang and Wang 2018b}
Huang J.W., Wang J.J. Higher order asymptotic behaviour of partial maxima of random sample from generalized Maxwell distribution under power normalization. Applied Mathematics-A Journal of Chinese Universities, 2018b,33(2): 177-187.

\bibitem{Resnick}
Resnick, S.I. (1987). {\it Extreme Value, Regular Variation, and Point Processes.} New York: Springer.

\end{thebibliography}
\end{document}